# Exploring localization in nonlinear oscillator systems through network-based predictions

Charlotte Geier[1, a)] and Norbert Hoffmann[1, 2]
[1)]*Hamburg University of Technology, Department of Mechanical Engineering, Hamburg, Germany*
[2)]*Imperial College London, Mechanical Engineering, London, UK*

(Dated: 6 April 2025)

Localized vibrations, arising from nonlinearities or symmetry breaking, pose a challenge in engineering, as the resulting high-amplitude vibrations may result in component failure due to fatigue. During operation, the emergence of localization is difficult to predict, partly because of changing parameters over the life cycle of a system. This work proposes a novel, network-based approach to detect an imminent localized vibration. Synthetic measurement data is used to generate a functional network, which captures the dynamic interplay of the machine parts, complementary to their geometric coupling. Analysis of these functional networks reveals an impending localized vibration and its location. The method is demonstrated using a model system for a bladed disk, a ring composed of coupled nonlinear Duffing oscillators. Results indicate that the proposed method is robust against small parameter uncertainties, added measurement noise, and the length of the measurement data samples. The source code for this work is available at[1].

**High-amplitude localized vibrations in engineering systems can lead to severe damage or failure, particularly in fields like aerospace and turbo-machinery. This study introduces a data-driven method that uses network-based techniques to anticipate these vibrations by analyzing system dynamics, offering insights that might help develop an early warning to prevent potential failures.**

---

## I. INTRODUCTION

Vibration localization poses an important challenge to the design and safe operation of engineering systems today. These high-amplitude vibrations of one or several parts of a machine may lead to fatigue, which may result in catastrophic system failure. This phenomenon is well known in turbo-machinery[2,3] and in the aerospace and aeronautical industries[4]. Due to the growing need for energy efficiency, machine parts are often designed to be more lightweight and slender structures[3]. As these systems are prone to exhibit vibration localization, this aspect remains an active field of research. Typical systems affected by localization include bladed disks[2,5], turbines[6–10], satellites, reflectors and antennas[4] and wind turbines[11]. Localized vibrations may occur in both linear and nonlinear systems. In linear systems, small inhomogeneities break the symmetry of a structure, causing localization, as observed by Anderson[12]. In engineering, imperfections can arise from manufacturing tolerances or wear during operation[11]. The effect is often called "mistuning" in the context of cyclic structures, such as turbo-machinery[2,5,13–17]. Localization can also arise in a perfectly symmetric structure if a nonlinearity is present in the system[18]. It has been shown that vibration localization can often be associated with the multi-stable behavior of a single oscillator in a symmetric chain[19,20]. Most real-world engineering systems include some form of nonlinearity, due to, for example, large deflections[3,11] or material nonlinearities[21,22].

Circumventing these high-amplitude vibrations plays a key role in engineering design[23]. Modeling and analyzing localization effects is often a non-trivial endeavor[24], as even small parameter variations may play a crucial role[11], and nonlinear effects are difficult to grasp. It is often difficult to assess the probability of a localized vibration occurring during operation, especially as an exact determination of the system parameters, which may change over a life cycle, is impossible. Therefore, an online early-warning system, which indicates an increasing risk of harmful vibrations, could be of interest to facilitate predictive maintenance[25,26] and avoid catastrophic events.

In this work, we develop a method to detect the emergence of a localized vibration due to parameter changes in an engineering system from measurement data. The approach is purely data-based and leverages methods founded in network science. Previous work has shown that different network-based approaches to time series analysis can track dynamical transitions, such as the emergence of synchronization scenarios[27–29] and dynamical regime shifts[30]. Recently, we proposed a method for inferring a functional network to analyze dynamical engineering systems[31]. A functional network encodes dynamical relationships between the elements of a system beyond their geometrical coupling. Properties of this functional network are leveraged to detect an imminent localized vibration. The method is presented using synthetic displacement data from a classical model for a bladed disk. The model consists of a set of cyclically coupled nonlinear Duffing oscillators. Results indicate that the procedure is robust against measurement noise, variability in the initial conditions, and parameter uncertainties.

This paper is structured in the following way: Section 2 introduces the model system and the data generation approach, as well as the network-based analyses. The results are shown in Section 3, along with studies regarding the robustness of the method. A discussion of the results follows in Section 4.

---

[a)]charlotte.geier@tuhh.de; http://www.tuhh.de/dyn

The work is completed by a conclusion in Section 5.

## II. METHODS

This work aims at detecting an imminent localized vibration from time series data from a dynamical system. An overview of the procedure is given in Fig. 1. In the first step, time series data is obtained from a model system, as illustrated in the following Subsection II A. The functional network is generated leveraging a recurrence-based approach, as proposed in[31] and presented in Subsection II B. The functional network is subsequently analyzed in terms of node indegrees and strongly connected components, as explained in Subsection II C. The procedure is repeated for a set of 100 random initial conditions. The full Python source code for this project is available at[1] under a GNU General Public License v3.0.

### A. Coupled Duffing oscillators

Coupled nonlinear Duffing oscillators form a model often used to represent turbine blades or bladed disks[11]. This dynamical system is known to exhibit rich dynamical behavior[32], which has been studied in detail for example in[32–34]. The model adopted in this work consists of $N = 10$ nonlinear Duffing oscillators with harmonic forcing. The oscillators are coupled to their next neighbors, forming a ring-like structure as shown in Fig. 1. The dynamics of the structure can be described by the second-order differential equation

$$\mathbf{M}\ddot{\mathbf{x}} + \mathbf{D}\dot{\mathbf{x}} + \mathbf{K}_l\mathbf{x} + \mathbf{F}_{nl} = \mathbf{f}(t), \quad (1)$$

where $\mathbf{M}, \mathbf{D}$, and $\mathbf{K}_l \in \mathbb{R}^{N \times N}$ are the mass, damping, and stiffness matrix, respectively. The displacement of each oscillator is given by $\mathbf{x} \in \mathbb{R}^{N \times 1}$, while $\dot{\mathbf{x}}$ and $\ddot{\mathbf{x}}$ denote the velocities and accelerations, respectively. The nonlinear terms are represented by $\mathbf{F}_{nl} \in \mathbb{R}^{N \times 1}$ and the forcing is given by $\mathbf{f}(t) \in \mathbb{R}^{N \times 1}$. In detail, the mass, damping, and stiffness matrices as well as the nonlinear terms and forcing are given in Appendix A. For the given set of parameters, every single oscillator exhibits bistable behavior and can oscillate with either a high or a low amplitude depending on the initial conditions, resulting in a multi-stable overall model. More details on the behavior of a single oscillator are given in[31]. In this study, one mass $m_4$ is assumed to decrease from 100 % to 80 %, e.g. due to material loss due to friction. The mass parameter is varied in 100 steps in $m_4 = [0.8, 1]$. As $m_4$ decreases, a localized vibration emerges at this oscillator. It is the aim of this work to detect the occurrence of this localized vibration before it occurs. This procedure could be useful in early warning methods or for predictive maintenance.

Time series measurements of the displacement of each oscillator are obtained from the model in Eq. 1 using an explicit Runge-Kutta method of order 4(5), which is implemented as 'Dopri5'[35] in Python Scipy[36]. The data samples of length 10 s are discretized at a step size $dt = 0.05$ s. For each variation of $m_4$, the measurement is obtained for two different sets of $M = 100$ initial conditions $\mathbf{x}_{0,1}$ and $\mathbf{x}_{0,2}$, to account for uncertain initial conditions that may be present in future real-word applications. The initial conditions are drawn from a random uniform distribution $\mathbf{x}_{0,1} \in [0, 0.1]$ and $\mathbf{x}_{0,2} \in [0, 0.01]$.

### B. Functional network

In engineering, interactions between machine parts are predominantly analyzed based on their geometric coupling properties. If two parts are physically connected, for example, as described by Eq. 1, their dynamics are assumed to be closely related. However, it has been observed that changes in one part of a machine, such as the tightening of a screw, may have significant effects on another part of the system further away. These relationships are named *functional relationships* for the sake of this work. Functional networks are designed to represent these relationships between system elements, complementary to their geometric coupling. Consequently, each node in the network represents one of the system machine parts or, in this case, one Duffing oscillator. The links between the nodes encode functional dependencies inferred from time series measurements.

The generation of the functional network from time series measurements follows the method introduced by the authors in[31]. This section will give a brief overview of the procedure. For more details on the procedure and its applicability, the interested reader may refer to[31].

First, a network with $N$ nodes and zero links is set up. The inference of functional coupling direction between the nodes is based on the idea of inter-system recurrence networks presented by Feldhoff et al. in[37]. To deduce the relationship between two nodes $i$ and $j$, the respective displacement measurements $x_i$ and $x_j$ are taken into account. Computation based on energy-related metrics would also be conceivable, but the use of oscillatory data $x$ enables a more direct application to real-world systems. The two time series are embedded into a high-dimensional space defined by their cross-recurrence. The cross-recurrence captures the mutual recurrence of the two time series in phase space. In the second step, the cross-recurrence matrix is interpreted as the adjacency matrix of a large network, allowing for the implementation of two network-based measures. The cross-clustering and cross-transitivity encode the number of neighbors who are also neighbors, and the number of cross-triangles over the number of cross-triples. Assuming a diffusive relationship between two elements, these measures can be interpreted as representing the functional relationship between two elements[37]. By obtaining the cross-recurrence and cross-transitivity from the high-dimensional space, coupling information is made available in the form of two scalar values. The sign of these measures indicates the coupling direction: $i \to j$, $i \leftarrow j$, or $i \leftrightarrow j$. The procedure is repeated for every possible combination of nodes, resulting in a directed, fully connected network. The computation of the functional networks is performed using the Python package pyunicorn[38].

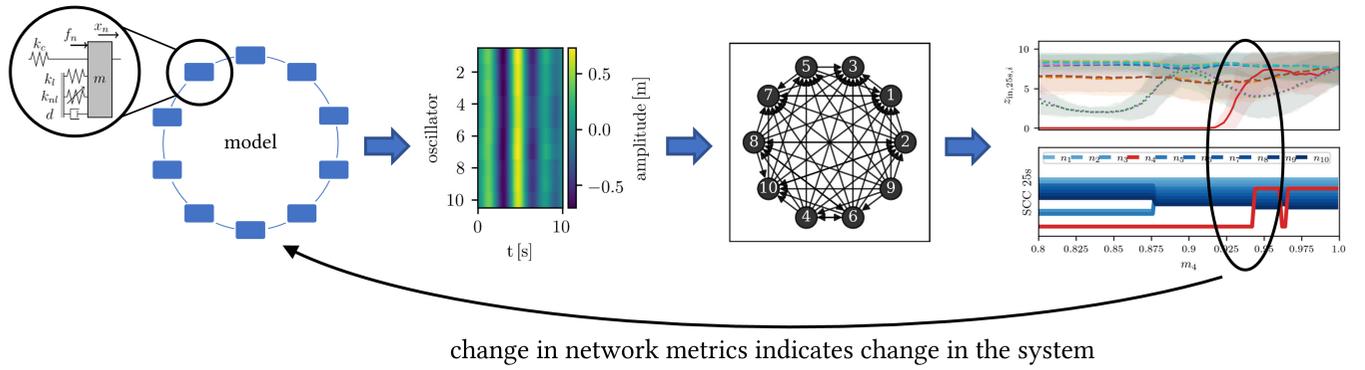

FIG. 1. Schematic of the analysis performed in this work. A symmetric model system with ten forced nonlinear Duffing oscillators with nearest-neighbor coupling is used to demonstrate the method. For this study, the mass of the fourth oscillator is decreased step-wise from 100% to 80%, introducing inhomogeneity into the system. Synthetic time series data is obtained from the model system in the form of 10 seconds of displacement measurements from each oscillator. A functional network is generated from the time series data using recurrence-based and network-based measures. The network encodes dynamical relationships between the oscillators beyond geometrical coupling. Each node represents an oscillator, while the links encode dynamical relationships between the nodes. The network is studied via the node in-degree, which describes the number of links incoming into each node. A statistical evaluation of the measure for 100 different initial conditions is performed. In addition, the formation of strongly connected components (SCCs) within the networks is studied. The evolution of these network measures enables the tracking of an imminent localized vibration

### C. Network analysis

The functional network is a directed, fully connected network with $N$ nodes. In this work, two indicators are used to analyze the network and infer information on the underlying dynamical system: the node-wise in-degree and the formation of strongly connected components (SCC). Both measures are described in this section.

The node-wise in-degree $z_{\text{in},i}$ counts the number of links that enter the $i$th node[39]. In this work, the evaluation is performed over $M = 100$ initial conditions per parameter variation. The mean value $z_{\text{in},i} = \frac{1}{M}\sum_M z_{\text{in},i,x_0}$ of the in-degrees are shown as a line, with the standard deviation $\sigma = \sqrt{\frac{1}{M}\sum_M (z_{\text{in},i,x_0} - z_{\text{in},i})^2}$ as a shaded area. While the node degree is an indicator of the importance of a single node in contributing to the system dynamics, it does not yield detailed information on the network structure. Therefore, the formation of SCCs is also taken into account in the analyses.

The evolution of SCCs within the network is detected using the "strongly connected components (SCC)"-algorithm introduced by Tarjan in[40] and further developed by Nuutila in[41]. The algorithm finds SCCs inside a network by combining sets of so-called path-equivalent nodes. Two nodes $i$ and $j$ are path equivalent, if there is a path from $i$ to $j$ and vice versa. For example, three nodes connected with all bi-directional edges would merge into a single SCC. The algorithm is available in the Python package networkX[42].

The following section is dedicated to representing the results obtained from these analyses.

### III. RESULTS

This section presents the evolution of node in-degree and SCCs as the mass parameter $m_4$ decreases. The studies are performed for two sets of initial conditions within different ranges. An analysis of the robustness against parameter uncertainties, noise, and sample time series length rounds off the studies.

The results obtained from time series data with initial conditions in range $\mathbf{x}_{0,1} \in [0, 0.1]$ are presented in Figure 2. The top row shows system dynamics for exemplary values of the mass parameter $m_4 = 0.8$ (left) $m_4 = 0.903$ (middle two panels), and $m_4 = 1.0$ (right panel). Each image represents 10 seconds of the dynamics of all ten oscillators, amplitudes are given in color-code. The means and standard deviation of each node in-degree $z_{\text{in},i}$, evaluated for all initial conditions, are presented in the second row. The third row displays the evolution of SCCs for one exemplary initial condition $\mathbf{x}_{0,a}$, see Appendix B. In this representation, each colored line represents a node. Touching lines represent the formation of a strongly connected component. For example, on the right-hand side of the figure, all nodes belong to the same strongly connected component, hence all lines lie close together. In the center of the figure, two SCCs have formed: node four separates from all other nodes. Exemplary networks representing specific SCCs are given in the bottom row.

For the uniform oscillation of the symmetric system with $m_4 = 1$ on the right-hand side of the figure, all nodes have the same mean in-degree $z_{\text{in},i}$, with a relatively large standard deviation. The corresponding network forms a single SCC, as indicated in the third row and the bottom image. As the mass of the fourth oscillator decreases, the in-degree of the respective node four also declines. While the first zero in-degrees for node four can be observed at $m_4 = 0.956$, the mean value $z_{\text{in},4}$ drops to zero at $m_4 = 0.931$. At the same time, the standard deviation for the mean in-degree of this node decreases. The first sustained localized vibration appears at $m_4 = 0.903$. The top panels illustrate the dependence of the dynamics on the initial conditions by showing two exemplary dynamics with different starting values. The left one localizes, but the right one doesn't. The network structure starts to disintegrate, reflecting this phenomenon. The node related to the localizing oscillator splits off from the rest of the network in most cases, while the remainder of the network still forms a variety of SCCs. One exemplary network structure is again shown in the bottom image. In the transition period from one large SSC to node four forming a separate component, node four switches between groups. This phenomenon reflects the use of transient dynamics in the computation, which leads to jumps in the distinction between localizing and non-localizing states when only a single initial condition is evaluated. Although this effect is shown only for one exemplary initial condition $\mathbf{x}_{0,a}$, it can be observed across the different initial conditions. Below this $m_4 = 0.877$, the localized vibration appears for every initial condition within the given range. The in-degree of the nodes corresponding to the neighboring oscillators in the model, nodes three and five, begin to drop. The network decomposes into three SCCs, the first formed by node four, the second by nodes three and five, and the third by the remainder of the network. The bottom image depicts the resulting structure. Note that the evolution of the SCC is shown only for a specific initial condition $\mathbf{x}_{0,a}$. For the remaining initial conditions, the evolution of component formation is not necessarily as distinct.

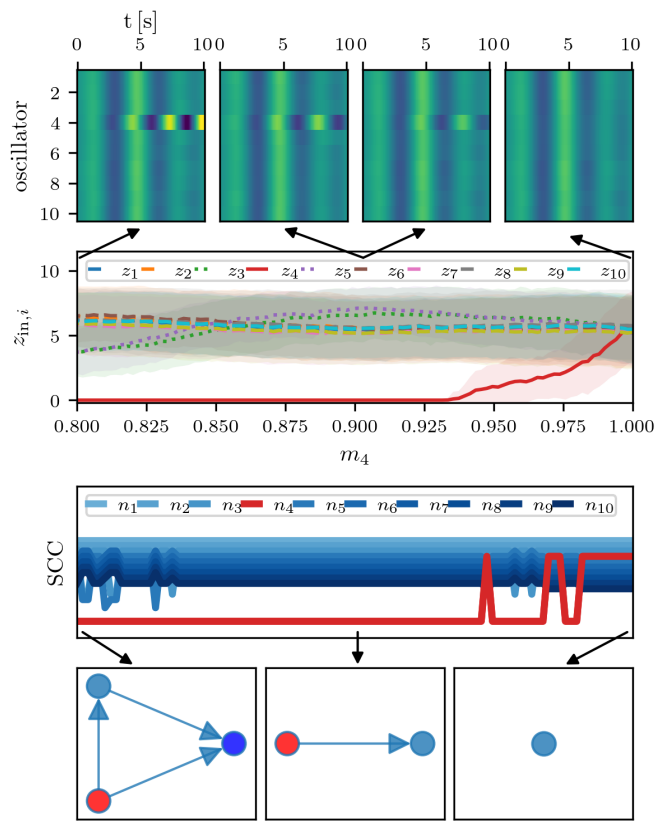

FIG. 2. Evolution of node in-degree and SCC formation for a decreasing mass $m_4$ of the fourth oscillator. The statistical evaluation is performed over a range of 100 initial conditions $\mathbf{x}_0 \in [0, 0.1]$. On the top, exemplary model time series taken at specific values of $m_4$ are presented. The color scale is identical to that in Fig. 1. For the given initial conditions, no localized vibrations occur for $m_4 > 0.903$. In a range of $0.877 < m_4 < 0.903$, localization can be observed sometimes, while for $m_4 < 0.877$, the localized vibration appears independent of the initial condition for the given range. The second row shows the evolution of the mean and standard deviation of the in-degree for each node. The mean in-degree of node four, $z_{\text{in},4}$ drops to zero at around $m_4 = 0.931$, while the first zero in-degrees for node four can be observed at $m_4 = 0.956$. The third row illustrates the evolution of SCCs for one illustrative initial condition $\mathbf{x}_{0,a}$, see Appendix B. The fourth node starts to split off the main network component at $m_4 = 0.966$, resulting in a one-node strongly connected component that singles out the mass that is affected by the localization. At $m_4 = 0.832$, nodes three and five start to split off as well. Small sub-figures below are exemplary networks condensed to their components, from each setting for comparison

Figure 3 presents the results for the second set of initial conditions $\mathbf{x}_{0,2} \in [0, 0.01]$. The given measures are analogous to the ones given in Figure 2. The results appear roughly similar to the ones observed for the larger range of initial conditions.

However, the dynamical transitions appear at different parameter values. For $0.880 < m_4$, no localized vibration can be observed. For $m_4 < 0.877$, the localized vibration emerges independent of the initial condition. The top panels of Figure 3 show the corresponding dynamics. The mean of the node-wise in-degree $z_{\text{in},i}$ exhibits a steeper flank compared to the image obtained with a larger range of initial conditions. The bands of standard deviations are much narrower. Compared to the larger range of initial conditions in Fig. 2, the in-degree $z_{\text{in},4}$ exhibits a transition period with a peak at $m_4 = 0.956$. This phenomenon is related in the increased sensitivity of the mean due to the smaller initial value range, which results in behavior similar to the evolution observed for in the SCC for a single initial condition. The in-degree of the fourth node, $z_{\text{in},4}$ drops to zero at $m_4 = 0.941$. Shortly before that, at $m_4 = 0.945$, the node splits off from the rest of the network to form its own strongly connected component. For an intermediate range of $0.822 < m_4 < 0.907$, nodes three and four also detach from the main SCC of the network to form separate or combined components. At the same time, the in-degree of these two nodes rises. The remainder of the parameter range exhibits a clear distinction between the giant SCC of the network and a small component formed by node four, which corresponds to the localizing oscillator. The evolution of network components is shown for one exemplary initial condition $\mathbf{x}_{0,b}$. In this case of initial conditions within a smaller range, the progression of SCC is much more consistent for different initial conditions, making it easier to extract a pattern.

Figure 4 illustrates the robustness of the approach to parameter uncertainties, measurement noise, and different time series lengths. The top panel of shows the node in-degrees obtained from a system with 1% random variation of the model

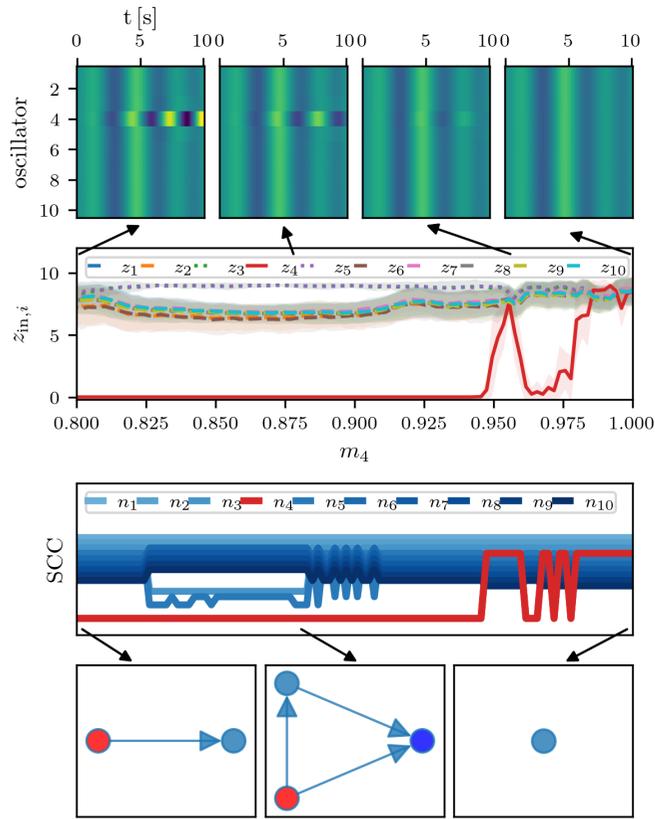

FIG. 3. Node in-degree and SCC over a variation of $m_4$ for the second set of initial conditions $x_0 \in [0, 0.01]$. The top row shows exemplary dynamics for distinct values of the parameter $m_4$ with a color scale identical to that in Fig. 1 Dynamic transitions are shifted to no localization for $m_4 > 0.88$ and always localization for $m_4 < 0.877$. The node in-degrees in the second row exhibit narrower standard deviation bands compared to the previous results. The evolution of SCC is given for $\mathbf{x}_{0,b}$ (see Appendix B). Initially, the uniform oscillation of the symmetric system results in identical mean in-degrees for all nodes and one SCC. The in-degree of node four, $z_{in,4}$, reaches zero at $m_4 = 0.941$ after a transition period similar to the evolution observed in the SSC, related to increased sensitivity due to the smaller initial value range. At $m_4 < 0.945$, the fourth node detaches from the main network SCC, as shown in the third. For further decrease of the parameter $m_4$, node four remains separated. For an intermediate period of $0.824 < m_4 < 0.907$, nodes three and five detach from the main network component to form separate or combined SCC. At the same time, the in-degrees of these two nodes are elevated. The bottom row shows networks illustrative of the three cases: a single SCC on the right, a three-component network in the middle, and a two-node condensed network on the left

parameters $m_i, k_i, k_{nl,i}$. The initial conditions $\mathbf{x}_0 \in [0, 0.1]$ are the same as in Figure 2. Compared to the results without parameter variations, the node in-degrees $z_{in,pu}$ diverge, even for an almost symmetric system where $m_4 = m_i = 1$. The standard deviation bands are broader, as expected from a more varied time series input. The primary indicator is still discernible: The node in-degree $z_{in,4}$ drops to zero around $m_4 = 0.927$, slightly later than in the previous result, where the node de-gree reached zero at $m_4 = 0.931$.

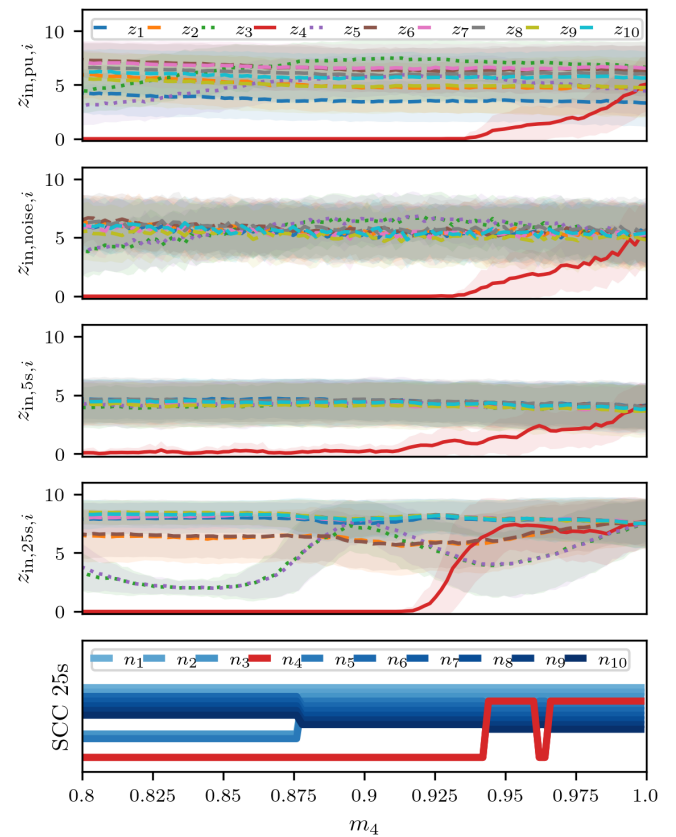

FIG. 4. Studies of the robustness of the approach. The top panel shows the resulting node in-degrees $z_{in,pu,i}$ when the underlying model system contains a 1% random variation of the system parameters $m_i, k_i, k_{nl,i}$. The initial conditions for this study are chosen such that $\mathbf{x}_0 \in [0, 0.1]$. While the standard deviation gets larger, and the individual node in-degree covers a larger range, the main observation of the degree of node four drops to zero is still possible. Results for noise resilience are shown in the second panel. Gaussian white noise with an amplitude of 5% of the standard deviation of the measurement time series is added to the displacement data. The resulting nodal in-degrees $z_{in,noise,i}$ are less smooth but still allow for the crucial analysis. The middle panel presents the node in-degrees $z_{in,5s,i}$ for a shorter time period. The bottom two panels depict the node in-degrees $z_{in,25s,i}$ computed and SCC evolution for 25 seconds input time series samples. The overall trend in both in-degree and SCC remains unchanged. The evolution of node in-degree of nodes three and five becomes more prominent, exhibiting a fluctuating behavior

The results for a noise-contaminated system are presented in the second panel of Figure 4. Noise is introduced into the system by adding Gaussian white noise with zero mean and a standard deviation of 5% of the standard deviation of the input time series to the time series data before generating the functional network. Excepting the noise, the parameter settings and initial conditions are kept constant from the first study in Figure 2. The resulting curves of mean in-degrees $z_{in,noise}$ are less smooth than the previous results. The standard deviation becomes slightly larger compared to the noise-free case.

The in-degree of node four, $z_{in,4}$ drops to zero at $m_4 = 0.923$, which indicates the imminent localized vibration as before.

The third and fourth panel show the resulting node in-degrees $z_{in,5s}$ and $z_{in,25s}$ if a time series of five and 25 seconds is used as input to the method, respectively. The bottom panel depicts the corresponding formation of SCCs for one exemplary initial condition $\mathbf{x}_{0,c}$ (see Appendix B) and time series length 25 seconds. The dynamical transitions are the same as illustrated in Figure 3. For $m_4 > 0.88$, no localized vibration exists within the given set of initial conditions. For $0.88 < m_4 < 0.877$, the localization appears in some cases, and for $m_4 < 0.877$, the localized vibration occurs for all given initial conditions. For the symmetric system with $m_4 = m_i = 1$, all nodes have the same mean in-degree, as is expected from previous results, and the network forms one single component. The mean in-degree $z_{in,5s,4}$ decreases less steeply than in the previous cases and never reaches zero, but drops below 0.1 at $m_4 = 0.874$, illustrating that the results, albeit less pronounced, remain the same. The mean node in-degree $z_{in,25s,4}$ and the network structure stays almost constant up to $m_4 = 0.95$. At $m_4 = 0.911$, the node in degree $z_{in,25s,4}$ drops to zero, while the standard deviation vanishes. At the same time, the in-degrees of the nodes corresponding to the neighboring oscillators three and five, $k_{in,25s,3}$ and $k_{in,25s,5}$ exhibit an oscillating path, dropping before $z_{in,25s,4}$ drops, then rising to a peak at $m_4 = 0.895$, then dropping again. The remaining in-degrees evolve pairwise according to their proximity to the localizing element. The evolution of SCCs shows the decomposition of the network into distinct parts evolves similarly as before: at $m_4 = 0.941$, node four splits off, and at $m_4 = 0.875$, the nodes related to the geometrically neighboring machine parts split off as well. The image resembles previous results but paints a picture with sharper transitions.

## IV. DISCUSSION

The results in the previous section indicate that the mean in-degree of nodes from a functional network can provide information on an imminent localized vibration for a system with one varying parameter. The decomposition of the network into distinct SCCs underlines this finding. Both results are valid for small parameter uncertainties and noise contamination of the measurement time series. The range of initial conditions underlying the time series data has an impact on the results, which become more distinct for a smaller range of initial conditions: The flanks of dropping in-degree become steeper, standard deviations shrink, and the SCC formation is more even. In a real-world application, this means that the range of initial conditions has to be bounded to some small value, such that each new measurement starts with similar starting values. A longer input time series results in a decreased standard deviation in node in-degrees and more distinct component formation compared to a shorter time series section. However, a relatively short section of a time series of 10 seconds, which covers about three periods of the system oscillation, is enough to detect the localized vibration. As the computation time increases with the number of time steps in each time series, using a shorter time series is computationally more efficient. However, it might also be possible to bootstrap the time series data, similar to approaches presented in[37], and only use a subset of the measurement points to decrease computation time. One advantage of this method is that it is entirely model-free. It requires only measurement data from the system in its "healthy" state and a regular repetition of these measurements. The data samples can be relatively short, but many samples may be necessary to smooth out the effects of uncertain initial conditions.

Although only results for a decreasing mass parameter were presented in this paper, the method works analogously for an increasing mass parameter. In this case, the in-degree of the node rises instead of drops, indicating the evolution of the parameter in the opposite direction. If the order of computation were reversed and the system were to evolve from a localized state towards a uniform oscillation, the same changes in degree and SCC would be observed, without hysteresis effects.

We expect our method to be transferable to slower, gradual localization phenomena which would manifest as more gradual transitions, and to multiple-oscillator localization, in which we anticipate a cluster of nodes to detach from the remaining network. The method does not currently distinguish between qualitative types of dynamics of each oscillator and would therefore likely not differentiate clearly between periodic localization and chimera states[31].

## V. CONCLUSION

In this work, we have presented a method for detecting an imminent localized vibration in a nonlinear cyclic oscillator chain, representing an engineering system model. This purely data-based approach does not depend on a physics-based model but requires measurement data from the system at hand. In the first step, a functional network is inferred from measurement data, which gives insight into the dynamic interplay of the system elements. In the second step, the functional network is analyzed, and the node in-degree is leveraged to detect the imminence and location of a localized vibration. The method is demonstrated using a cyclic chain of nonlinear Duffing oscillators.

We have shown that the approach can be applied to detect the onset of a localized vibration even before it is discernible in the time series data. SCCs in the functional network can be related to dynamical relationships between the machine parts, which can often be related to the symmetries within the vibration. As a localized vibration emerges, the node related to the part that oscillates at high amplitude separates from the remaining network, forming its separate component. This observation enables us to determine the location of the localized vibration. A statistical evaluation smooths over the effect of uncertain initial conditions. The picture becomes sharper if the initial conditions lie within a smaller range. At the same time, the SCCs become more distinct, and the evolution of the network SCCs over a parameter change can be discerned more clearly. The results indicate that the method works well even under uncertain conditions, such as unknown initial con-ditions, small parameter uncertainties, and the presence of noise.

Further studies will include the detection variations in other model parameters, such as the coupling stiffness, and the application to larger systems with a more complex coupling structure. We hope this work might inspire other work on leveraging functional relationships for analyzing dynamical systems.

## ACKNOWLEDGMENTS

C.G. is thankful to the DFG (German Research Foundation) for support through project number 510246309.

## DATA AVAILABILITY STATEMENT

The data and code that support the findings of this study are openly available in Zenodo at https://doi.org/10.5281/zenodo.12611988, reference number 12611988.

## CONFLICT OF INTEREST STATEMENT

The authors have no conflicts to disclose.

## Appendix A: Model details

In detail, the mass, damping, and stiffness matrices as well as the nonlinear terms and forcing of Equation 1 are given by

$$\mathbf{M} = \begin{bmatrix} m_1 & & \\ & \ddots & \\ & & m_{10} \end{bmatrix}, \mathbf{D} = \begin{bmatrix} d_1 & & \\ & \ddots & \\ & & d_{10} \end{bmatrix},$$

$$\mathbf{K}_l = \begin{bmatrix} 2k_l & -k_c & & & -k_c \\ -k_c & 2k_l & -k_c & & \\ & & \ddots & & \\ & & -k_c & 2k_l & -k_c \\ -k_c & & & k_c & 2k_l \end{bmatrix}, \quad (A1)$$

where $m_i = 1$ is the mass of each oscillator, $d_i = \alpha m_i$, $\alpha = 0.1$ describes the damping, $k_l = 1$ the linear spring stiffness and $k_c = 0.1$ the coupling stiffness terms. The nonlinearity is introduced as a cubic stiffness and the forcing is given as harmonic forcing

$$\mathbf{F}_{nl} = \begin{bmatrix} k_{nl} x_1^3 \\ \vdots \\ k_{nl} x_{10}^3 \end{bmatrix}, \mathbf{f}(t) = \begin{bmatrix} F\cos(\Omega t) \\ \vdots \\ F\cos(\Omega t) \end{bmatrix}, \quad (A2)$$

with nonlinear spring stiffness $k_{nl} = 2$, forcing amplitude $F = 1$ and forcing frequency $\Omega = 2$.

## Appendix B: Initial conditions for numerical data

Initial conditions for the underlying time series data used to showcase the evolution of SCCs in Figure 2

$x_{0,a} = [0.0678, 0.0392, 0.0330, 0.03074, 0.0738, 0.0672, 0.0413,$
$\quad 0.0780, 0.06037291, 0.0747, 0.0062, 0.0381, 0.0933,$
$\quad 0.0543, 0.0505, 0.0356, 0.0193, 0.0180, 0.0577, 0.0390],$

in Figure 3

$x_{0,b} = [0.0092, 0.0089, 0.0045, 0.0023, 0.0013, 0.0018, 0.008,$
$\quad 0.0065, 0.0002, 0.0048, 0.0026, 0.005, 0.0087, 0.0033,$
$\quad 0.0006, 0.0019, 0.0041, 0.002, 0.005, 0.0007],$

and in Figure 4d

$x_{0,c} = [0.0466, 0.0039, 0.0121, 0.0876, 0.0077, 0.0366, 0.0205,$
$\quad 0.0105, 0.0442, 0.0656, 0.0261, 0.0725, 0.0919, 0.0103,$
$\quad 0.0613, 0.0310, 0.0878, 0.0022, 0.0736, 0.0126].$